\documentclass[letterpaper, 10 pt, conference]{ieeeconf}  
\IEEEoverridecommandlockouts                             
\usepackage{graphics} 
\usepackage{epsfig} 
\usepackage{times} 
\usepackage{amsmath} 
\usepackage{amssymb} 
\usepackage{color}
\usepackage{algorithm}
\usepackage{algpseudocode}
\usepackage{bbm}
 \usepackage{url}
 \usepackage{bbold}
\usepackage{mathtools}
\usepackage[table]{xcolor}
\usepackage{shorthands}
\newcommand{\mc}{\mathbb}
\usepackage{subcaption}
\usepackage{flushend}
\usepackage{cite} 

\definecolor{babyblueeyes}{rgb}{0.63, 0.79, 0.95}
\definecolor{bluegray}{rgb}{0.4, 0.6, 0.8}
\definecolor{aliceblue}{rgb}{0.94, 0.97, 1.0}
\definecolor{carolinablue}{rgb}{0.6, 0.73, 0.89}
\definecolor{columbiablue}{rgb}{0.61, 0.87, 1.0}
\definecolor{cyan(process)}{rgb}{0.0, 0.72, 0.92}
\definecolor{lavenderblue}{rgb}{0.25, 0.7, 0.7}
\title{\LARGE \bf
Oracle Complexity Reduction for Model-free LQR:\\ A Stochastic Variance-Reduced Policy Gradient Approach
}

\author{Leonardo F. Toso, Han Wang, and James Anderson
\thanks{This material is based upon work supported in part by  NSF awards 2144634 \& 2231350. Leonardo F. Toso is funded by the Columbia Presidential Fellowship. The authors are with the Department of Electrical Engineering,  Columbia University in the City of New York, New York, NY, 10027, USA. Email: \texttt{\{lt2879, hw2786, james.anderson\}@columbia.edu}.} 
\thanks{
        {\tt\small }}%
\thanks{
        {\tt\small }}%
}

\begin{document}

\allowdisplaybreaks
\maketitle
\thispagestyle{empty}
\pagestyle{empty}

\begin{abstract}
We investigate the problem of learning an $\epsilon$-approximate solution for the discrete-time Linear Quadratic Regulator (LQR) problem via a Stochastic Variance-Reduced Policy Gradient (SVRPG) approach. Whilst policy gradient methods have proven to converge linearly to the optimal solution of the model-free LQR problem, the substantial requirement for two-point cost queries in gradient estimations may be intractable, particularly in applications where obtaining cost function evaluations at two distinct control input configurations is exceptionally costly. To this end, we propose an oracle-efficient approach. Our method combines both one-point and two-point estimations in a dual-loop variance-reduced algorithm. It achieves an approximate optimal solution with only $\mathcal{O}\left(\log\left(1/\epsilon\right)^{\beta}\right)$ two-point cost information for $\beta \in (0,1)$.
\end{abstract}


\section{Introduction}

Policy gradient (PG) methods have attracted significant attention in model-free reinforcement learning (RL), in large part due to their simplicity of implementation. Within the context of control, and the LQR problem specifically (where analytic solutions are known), a lot of recent work has focused on connecting system theoretic properties such as controllability, with learning theoretic measures such as sample complexity~\cite{ziemann2022policy}. As first shown in~\cite{fazel2018global} and further analyzed in~\cite{gravell2020learning,malik2019derivative, mohammadi2020linear}, PG methods converge to the global optimal solutions despite the lack of convexity in the LQR problem. This significant result, combined with the adaptability of PG in the model-free setting, has opened up a line of research that addresses classical control problems using PG-based approaches \cite{hu2023toward, perdomo2021stabilizing}.

In the model-free LQR setting, policy gradient descent relies on a finite-sample estimate of the true gradient, often acquired through derivative-free (otherwise known as zeroth-order) methods. We refer the reader to ~\cite{malik2019derivative} for specific application of zeroth-order methods to LQR control and \cite{spall2005introduction} for general background. Zeroth-order gradient estimation approaches are particularly valuable for applications where the computational resources needed for exact gradient evaluations may be impractical, or when cost-query information is \emph{only} accessible through a black-box procedure. 

Despite providing flexibility by avoiding the explicit computation of gradients, zeroth-order gradient estimations with one-point (ZO1P) or two-point (ZO2P) queries frequently produce biased estimations accompanied by large variances \cite{malik2019derivative}. In order to counteract this, large sample sizes are required to accurately estimate the gradients.

Whilst ZO2P provides a reduced variance relative to ZO1P, it necessitates querying the cost function at two distinct control input configurations, which can be prohibitively impractical for certain applications (e.g., robot path planning \cite{peters2008reinforcement}). Addressing this limitation is crucial for developing efficient approaches applicable to real-world scenarios.  

Motivated by these challenges, one line of work focuses on leveraging data from multiple similar systems to mitigate variance and thereby reduce the sample complexity of policy gradient methods  \cite{wang202model, tang2023zeroth}. However, for the single-agent setting  it remains unclear how we can devise a more computationally efficient approach without resorting to second-order techniques.

On the other hand,  in supervised learning  and RL, SVRPG approaches have demonstrated their effectiveness in significantly reducing variance and enhancing sample efficiency for PG methods \cite{johnson2013accelerating, papini2018stochastic}. Such methods leverage the well-known control variate analysis, which incorporates both current and past gradient information to form a descent direction that reduces the estimation's variance. This concept motivates the following question addressed in this work:

\emph{Can we design an oracle-efficient solution for addressing the model-free LQR problem by building upon the success of stochastic variance-reduced approaches?}.

\noindent \textbf{Our Contributions}: Toward this end, our main contributions are summarized as follows: 

\begin{itemize}
    \item This is the first work to propose a stochastic variance-reduced policy gradient algorithm featuring a mixed zeroth-order gradient estimation scheme for tackling the model-free and discrete-time LQR problem.
    \item Theoretical guarantees demonstrate the convergence (Theorem \ref{theorem:convergence}) of our approach, while ensuring stability of the system under the iterated policy (Theorem \ref{theorem:stability}).
    \item We establish conditions on the problem parameters under which our approach achieves an $\epsilon$-approximate solution with $\mathcal{O}\left(\log\left(1/\epsilon\right)^{3-2\beta}\right)$ queries, while utilizing only $\mathcal{O}\left(\log\left(1/\epsilon\right)^{\beta}\right)$ two-point query information for $\beta \in (0,1)$. This oracle complexity improves upon the best known result $\mathcal{O}\left(\log\left(1/\epsilon\right)\right)$ by a factor of $\mathcal{O}\left(\log\left(1/\epsilon\right)^{1-\beta}\right)$ (Corollary \ref{corollary:oracle_reduction}).
\end{itemize} 

\noindent \textbf{Main result overview:} The SVRPG method we propose requires a slightly larger number of queries, specifically we require $\mathcal{O}\left(\log\left(1/\epsilon\right)^{3-2\beta}\right)$, (this includes one \emph{and} two-point queries) in comparison to $\mathcal{O}\left(\log(1/\epsilon)\right)$ required by the standard ZO2P approach, in order to achieve an $\epsilon$-approximate solution -- the difference is only  a logarithmic factor, for \emph{large} $\beta$. However, our approach requires considerably fewer two-point queries, specifically a factor of $\mathcal{O}\left(\log\left(1/\epsilon\right)^{1-\beta}\right)$ fewer, for \emph{small} $\beta$. This underscores the benefit of our technique, particularly in applications where conducting two-point function evaluations is prohibitively costly.

\subsection{Related Work}

\noindent \textbf{Model-free LQR via Policy Gradient Methods:} PG methods have been extensively explored as a solution to solve the model-free LQR problem in both discrete \cite{fazel2018global, gravell2020learning, malik2019derivative, mohammadi2020linear, ziemann2022policy, hu2023toward} and continuous-time settings \cite{mohammadi2019global, mohammadi2020random, mohammadi2021convergence}. Despite of the non-convexity of the LQR landscape under the policy search, Fazel et al. \cite{fazel2018global}  proved theoretical guarantees for the global convergence of PG methods for both model-based and model-free settings. Table \ref{table:sample_complexity_LQR} summarizes the sample complexity of the aforementioned work.

Although there has been an evident sample complexity reduction from $\mathcal{O}(\frac{1}{\epsilon} \log\left(1/\epsilon\right))$ \cite{malik2019derivative} to $\mathcal{O}\left(\log\left(1/\epsilon\right)\right)$ \cite{mohammadi2020linear}, this is primarily a result of a more refined analysis rather than algorithmic development.\footnote{We use big-O notation $\mathcal{O}(\cdot)$ to omit constant factors in the argument.} In this work, we propose a SVRPG algorithm to reduce the number of two-point queries required to obtain an $\epsilon$-approximate solution for the LQR problem.\vspace{0.2cm}

\noindent \textbf{Stochastic Variance-Reduced Policy Gradient:} Stochastic variance-reduced gradient descent (SVRG) have emerged as a sample-efficient  solution technique for  non-convex finite-sum optimization problems. Whilst SVRG methods have long been established for non-convex optimization problems (e.g., SVRG \cite{johnson2013accelerating}, SAG \cite{roux2012stochastic}, and SAGA \cite{defazio2014saga}), their extension to online RL settings is a relatively recent development (e.g., SVRPG \cite{papini2018stochastic, xu2020improved, liu2020improved}). This extension has presented unique challenges, primarily stemming from policy non-stationarity and approximations in the computation of the gradient. Furthermore, SVRPG approaches generally rely on the assumption of unbiased gradient estimation, a condition that rarely holds for derivative-free techniques. This has been addressed in \cite{ji2019improved, liu2018zeroth} for finite-sum, non-convex problems.

We emphasize that our work does not revolve around a simple extension of the results in \cite{papini2018stochastic, xu2020improved, liu2020improved} (online RL setting) or \cite{ji2019improved, liu2018zeroth} (non-convex finite-sum problem). In contrast to the latter, our LQR setting encompasses an online optimization problem with a single cost function. As a result, the sampling variance reduction benefit of using zeroth-order variance-reduced methods cannot be simply extended to our setting. On the other hand, in our setting we have the stabilizing policy requirement which is commonly taken for granted in the Markov Decision Process (MDP) case \cite{papini2018stochastic, xu2020improved, liu2020improved} with irreducibly and aperiodicity assumptions on the policy. Moreover, the zeroth-order gradient estimation produces biased estimations. This necessitates further derivations to control this bias as we will discuss later.

\begin{center}
\begin{table*}
\centering
\normalsize
\caption{Comparison on the sample complexity $(\mathbb{S}_c)$, and two-point oracle complexity ($\mathcal{N}_{\text{ZO2P}}$) required to achieve $\mc{E}\left(C(K_{\text{out}}) - C(K^*)\right) \leq \epsilon$. Here $\beta \in (0,1)$.}
\begin{tabular}{
c 
c  c }
\hline
\textbf{Methods}        & \textbf{Sample Complexity} 
 (\textbf{$\mathbb{S}_c$})   & \textbf{Two-point Oracle Complexity ($\mathcal{N}_{\text{ZO2P}}$)}   \\ \hline
 PG - ZO1P (Fazel et al (2018), \cite{fazel2018global}) &  $\mathcal{O}(1/\epsilon^{4}\cdot\log\left(1/\epsilon\right))$  &  -     \\ \hline
PG - ZO1P (Gravell et al (2019), \cite{gravell2020learning}) &   $\mathcal{O}(1/\epsilon^{4}\cdot\log\left(1/\epsilon\right))$  &  -     \\ \hline
PG - ZO1P (Malik et al. (2019), \cite{malik2019derivative}) & $\mathcal{O}(1/\epsilon^{2}\cdot\log\left(1/\epsilon\right))$  &  - \\ \hline
PG - ZO2P (Malik et al. (2019), \cite{malik2019derivative}) & $\mathcal{O}(1/\epsilon\cdot\log\left(1/\epsilon\right))$  &   $\mathcal{O}(1/\epsilon\cdot\log\left(1/\epsilon\right))$ \\ \hline
PG - ZO2P (Mohammadi et al. (2020), \cite{mohammadi2020linear}) & $\mathcal{O}(\log\left(1/\epsilon\right))$  & $\mathcal{O}(\log\left(1/\epsilon\right))$  \\ \hline
\rowcolor{aliceblue}
SVRPG - Algorithm \ref{algorithm:LQR_SVRPG} (This paper)  & $\mathcal{O}\left(\log\left(1/\epsilon\right)^{3-2\beta}\right)$  & $\mathcal{O}\left(\log\left(1/\epsilon\right)^{\beta}\right)$ \\ \hline
\end{tabular}
\label{table:sample_complexity_LQR}
\end{table*}
\end{center}

\vspace{-1.11cm}
\section{Preliminaries}
We  summarize key policy gradient results for the LQR problem as well as derivative-free optimization techniques.

\subsection{Discrete-time Linear Quadratic Regulator}

Consider the discrete-time LTI system 
\begin{align}\label{eq:LTI_system}
    x_{\tau+1}=Ax_{\tau} + Bu_{\tau},\;\ x_0  \stackrel{\text{i.i.d.}}{\sim} \mathcal{X}_0,
\end{align}
where $x_\tau \in \mathbb{R}^{n_x}$, $u_\tau \in \mathbb{R}^{n_u}$, and  $x_0$ denote the state and input at time $\tau$, and the initial condition. The optimal LQR policy associated with \eqref{eq:LTI_system} is $u_\tau=-K^*x_\tau$ where $K^*$ solves
\begin{align}\label{eq:LQR_problem}
 & \argmin_{K \in \mathcal{K}} \left\{ C(K) := \mc{E}_{x_0 \sim \mathcal{X}_0} \left[\sum_{\tau=0}^{\infty} x^{\top}_\tau Qx_\tau + u^{\top}_\tau Ru_\tau\right] \right\},\notag\\
& \hspace{1cm} \text{subject to}\:\ \eqref{eq:LTI_system}
\end{align}
where $Q \in \mathbb{S}^{n_x}_{\succ 0}$, $R \in \mathbb{S}^{n_u}_{\succ 0}$, and $\mathcal{K} := \{K| \rho(A - BK) < 1\}$ denotes all stabilizing  controllers $K \in \mathbb{R}^{n_u\times n_x}$. 
The optimal cost is assumed to be finite. This is satisfied when $(A, B)$ is controllable.

In the model-based setting  the optimal controller is given by:
$ K^* := -\left(R + B^\top P B\right)^{-1} B^\top P A, $
where $P \in \mathbb{S}^{n_x}_{\succ 0}$ is the solution of the
Algebraic Riccati Equation (ARE)~\cite{hewer1971iterative}. In the absence of the system model $(A,B)$, there is no way to implement an ARE-derived controller. Notably, motivated by the fact that traditional RL techniques aim to find optimal policies for unknown MDPs through direct exploration of the policy space, the line of work led by Fazel et al.  \cite{fazel2018global} and followed by \cite{gravell2020learning, malik2019derivative, mohammadi2020linear, mohammadi2020learning, mohammadi2021convergence, mohammadi2019global} have proved guarantees for the global convergence of PG methods for both model-based and model-free LQR. This is achieved by leveraging fundamental properties of the LQR cost function. Next, we revisit the updating rule of the model-free LQR problem through policy gradient, as well as its important properties. 

Suppose that instead of having the true gradient $\nabla C(K_l)$ at the $l$-th iteration, we posses a finite-sample estimate $\widehat{\nabla}C(K_l)$. The policy gradient method's update rule for the LQR problem can be expressed as follows:
\begin{align}\label{eq:updating_rule}
    K_{l+1} = K_l - \eta \widehat{\nabla}C(K_l), \quad l = 0,1,\ldots,L-1
\end{align}
 where $\eta$ represents a positive scalar step-size. We require the following standard assumption\cite{fazel2018global, gravell2020learning, malik2019derivative, mohammadi2020linear}.  
\begin{assumption} \label{assumption:stabilizing_K0} We have access to an initial stabilizing controller $K_0$ such that $\rho(A - BK_0) < 1$.
\end{assumption}

\begin{remark}
Note that if the initial controller $K_0$ fails to stabilize system \eqref{eq:LTI_system}, the PG in \eqref{eq:updating_rule} cannot iteratively converge to a stabilizing policy since $\widehat{\nabla} C(K_0)$ becomes undefined.
\end{remark}

\begin{definition}  The sublevel set of stabilizing feedback controllers 
$\mathcal{G} \subseteq \mathcal{K}$ is defined as follows
\begin{align*}
   \mathcal{G} := \{K\;\ |\;\ C(K) - C(K^*) \leq \xi \Delta_0\},
\end{align*}
where $\Delta_0=C(K_0) - C(K^*)$ and $\xi$ is any positive constant.
\end{definition}


\begin{lemma}  \label{lemma:Lipschitz_true_cost} Given two stabilizing policies $K^{\prime}$, $K \in \mathcal{G}$ such that $\|K^{\prime} -K\|_F  \leq h_{\Delta}(K)\ <\infty$, it holds that
\begin{align*}
&\left|C\left(K^{\prime}\right)-C(K)\right| \leq h_{\text {cost}}(K) C(K)\|K^{\prime} -K\|_F, \notag\\
&\left\|\nabla C\left(K^{\prime}\right)-\nabla C(K)\right\|_F \leq h_{\text {grad}}(K)\|K^{\prime} -K\|_F.  
\end{align*}
\end{lemma}

\begin{lemma}  \label{lemma:PL_condition}
Let $K^* \in \mathcal{G}$ be the optimal policy that solves \eqref{eq:LQR_problem}. Thus, it holds that
\begin{align*}
 C(K)-C\left(K^*\right) \leq \frac{1}{\lambda}\|\nabla C(K)\|_F^2, 
\end{align*}
for any stabilizing controller $K \in \mathcal{G}$.  
\end{lemma}

A detailed proof of the above lemmas, along with the explicit expressions for $h_{\Delta}(K)$, $h_{\text{cost}}(K)$, $h_{\text{grad}}(K)$, and $\lambda$, can be found in \cite{gravell2020learning}. We direct the reader to Appendix \ref{appendix:auxiliary_def} for the definition of $\bar{h}_{\text{grad}}$, $\bar{h}_{\text{cost}}$, and $\underline{h}_{\Delta}$ that are positive coefficients we use further in our derivations.

\subsection{Zeroth-Order Gradient Estimation}

Given a positive scalar smoothing radius, denoted as $r$, and randomly sampled matrices $U_1,\ldots, U_{m}$ drawn i.i.d. from the uniform distribution $\mathcal{S}_r$ of matrices with $\|U\|_F=r$, and considering a given stabilizing policy $K \in \mathcal{G}$, we define the one-point and two-point zeroth-order gradient estimations of the true gradient $\nabla C(K)$ as follows:

\begin{align*}
    \textbf{ZO1P}: \overline{\nabla}C(K) := \sum_{i=1}^{m} \frac{d C(K+U_i)U_i}{mr^2},
\end{align*}
\vspace{-0.2cm}
\begin{align*}
    &\textbf{ZO2P}:  \widetilde{\nabla}C(K) :=\sum_{i=1}^{m} \frac{d\left(C(K+U_i)-C(K-U_i)\right)U_i}{2mr^2},
\end{align*}
\noindent where $d=n_xn_u$ and $C(\cdot)$ denotes the true cost value provided by an oracle.

We emphasize that, in practice, we have a finite number of samples denoted by $m$ to compute ZO1P and ZO2P. Consequently, both ZO1P and ZO2P gradient estimation schemes exhibit an inherent bias. In addition, for  simplicity we assume access to the true cost, as provided by an oracle~\cite{malik2019derivative}. In reality, practical limitations prevent us from simulating our system over an infinite horizon. However, as in \cite[Appendix B]{gravell2020learning} the finite horizon approximation for the cost is upper-bounded by the true cost, with the approximation error controllable by the horizon length. Our  work can thus be readily extended to this finite-horizon approximated cost setting.

Moreover, the expressions of ZO1P and ZO2P shed light on the fact that whilst ZO2P requires more computational resources due to the need for two cost-query information for each sampled matrix $U \stackrel{\text{i.i.d.}}{\sim} \mathcal{S}_r$, it offers a lower-variance estimation, which results in a more efficient sample complexity, compared to ZO1P~ \cite{malik2019derivative}. This makes ZO2P a more favorable choice over ZO1P gradient estimation. Next, we present the PG algorithm with ZO2P gradient estimations for solving the model-free LQR.

\begin{algorithm}[H]
\caption{PG with ZO2P Gradient Estimation.} 
\label{algorithm:LQR_ZO2P}
\begin{algorithmic}[1]
\State \text {\textbf{Input:} } $L$,  $\eta$, $n_1$, $r$, $K_0$\
\State  \textbf{for} $l=0, \ldots, L-1$ \textbf{ do }\
\State \quad Compute $\widetilde{\nabla}C(K_l)$ with $r$ via ZO2P \
\State \quad $K_{l+1}=K_l-\eta \widetilde{\nabla}C(K_l)$ \
\noindent \State \textbf{ end for}\\
\textbf{Output} $K_{\text{out}}:=K_{L}$
\end{algorithmic}
\end{algorithm}

It is well-established~\cite{mohammadi2020linear} that under certain conditions on the quality of the estimated gradient, i.e., with $n_1$ large and $r$ small, Algorithm \ref{algorithm:LQR_ZO2P} converges linearly to the optimal solution of \eqref{eq:LQR_problem} while ensuring $K_l \in \mathcal{G}$ at each iteration. However, due to the still high variance of the gradient estimation step, the required number of two-point queries to achieve an $\epsilon$-approximate solution may become prohibitively large.

\vspace{-0.2cm}
\section{An SVRPG Algorithm for model-free LQR}

With the purpose of reducing the number of two-cost query information to achieve an $\epsilon$-approximate solution  we propose a SVRPG approach featuring a mixed gradient estimation scheme. The idea is to use  a ZO2P gradient estimate in the outer-loop and a ZO1P estimate in the inner-loop so as to lower the computational complexity associated with two-point cost queries compared to Algorithm \ref{algorithm:LQR_ZO2P}. The need for two-point cost query information arises only periodically instead of at each iteration. 

\begin{algorithm}
\caption{LQR via SVRPG} 
\label{algorithm:LQR_SVRPG}
\begin{algorithmic}[1]
\State  \textbf{Input:} $N$, $T$, $\eta$, $n_1$, $n_2$, $r_{\text{out}}$, $r_{\text{in}}$, $K_T^0:=\widetilde{K}^0:=K_0$.\
\State  \textbf{for} $n=0, \ldots, N-1$ \textbf{ do }\
\State \quad   $K_0^{n+1}:=\widetilde{K}^n:=K_T^n$ \
\State \quad Compute $\tilde{\mu} = \widetilde{\nabla}C(\tilde{K}^n)$ with $r_{\text{out}}$ \Comment{ZO2P} \
\State \quad  \textbf{for} $t=0, \ldots, T-1$ \textbf{ do } \
\State \quad \quad  Compute $\overline{\nabla}C(K^{n+1}_t)$, $\overline{\nabla}C(\tilde{K}^n)$ with $r_{\text{in}}$ \Comment{ZO1P} \
\State \quad \quad  $v_t^{n+1}=\tilde{\mu} + \overline{\nabla}C(K^{n+1}_t) - \overline{\nabla}C(\tilde{K}^n)$ \
\State \quad \quad $K_{t+1}^{n+1}=K_t^{n+1}-\eta v_t^{n+1}$ \
\State \quad  \textbf{ end for } \
\noindent \State \textbf{ end for}\\
\textbf{Output} $K_{\text{out}}:=K_{T}^{N}$.
\end{algorithmic}
\end{algorithm}

In contrast to Algorithm \ref{algorithm:LQR_ZO2P} our SVRPG algorithm  divides the total number of iterations into $N$ epochs, each of length $T$. For each epoch (outer-loop), we estimate gradients using $n_1$ samples with smoothing radius $r_{\text{out}}$, whereas inside each epoch (inner-loop) we use $n_2$ samples with smoothing radius $r_{\text{in}}$. In line 3, we fix the current policy $\Tilde{K}^n$ and compute $\widetilde{\nabla}C(\tilde{K}^n)$ via ZO2P. Throughout the inner-loop iterations, we estimate $\overline{\nabla}C(K^{n+1}_t)$ and $\overline{\nabla}C(\tilde{K}^n)$ with the same set of samples via ZO1P. Finally, in line 8 we perform a gradient descent step, using the stochastic variance-reduced gradient computed in line 7.

To close this section, we briefly discuss the idea behind SVRG-based methods. Consider a fixed stabilizing policy $\tilde{K} \in \mathcal{G}$ and estimate $\widetilde{\nabla}C(\Tilde{K})$ using $n_1$ samples. Then perform $K \leftarrow K - \eta v$, with
\begin{align*}
&v = \widetilde{\nabla}C(\Tilde{K}) + \overline{\nabla}C(K)- \overline{\nabla}C(\tilde{K}),
\end{align*}
where $\overline{\nabla}C(K)$ and $\overline{\nabla}C(\tilde{K})$ are estimated by using the \emph{same} set of $n_2$ samples. Note that $\mc{E}\widetilde{\nabla}C(\tilde{K}) = \mc{E}\overline{\nabla}C(\tilde{K})$ (see Appendix \ref{appendix:proof_bias}). Therefore, since $\overline{\nabla}C(K)$, and $ \overline{\nabla}C(\tilde{K})$ are correlated through their samples, the variance of the stochastic gradient $v_l$ might be reduced by controlling the covariance across the gradient estimations. That is, $\var(v) = \var(X-Y) = \var(X) + \var(Y) -2\cov(X,Y)$, with $X=\overline{\nabla} C(K)$, $Y=\widetilde{\nabla} C(\tilde{K}) - \overline{\nabla} C(\tilde{K})$,  and $\cov(\cdot,\cdot)$ denotes the covariance operator.

\section{Theoretical Guarantees}\label{sec:theoretical_guarantees}

 Without loss of generality and for the purpose of the theoretical analysis only,  set $r_{\text{out}} = r_{\text{in}} = r$ in Algorithm \ref{algorithm:LQR_SVRPG}. In  Proposition \ref{prop:linear_ZO2P} we first establish the convergence rate of  Algorithm \ref{algorithm:LQR_ZO2P}.  This allows for a fair comparison on the sample and oracle complexities of Algorithm 2, detailed in Corollaries \ref{corollary:sample_complexity} and \ref{corollary:oracle_reduction}. Moreover, we outline the conditions under which Algorithm \ref{algorithm:LQR_SVRPG} converge to the optimal solution (Theorem \ref{theorem:convergence}), all while staying within the stabilizing sub-level set (Theorem \ref{theorem:stability}) throughout the algorithm's iterations.

\begin{proposition} \label{prop:linear_ZO2P} (Convergence of Algorithm \ref{algorithm:LQR_ZO2P}) Suppose the smoothing radius, number of samples, and number of iterations are in the order of $n_1={\mathcal{O}}(1)$, $r=\mathcal{O}(\sqrt{\epsilon})$ and $L = \mathcal{O}\left(\log(1/\epsilon)\right)$, respectively. Then, Algorithm \ref{algorithm:LQR_ZO2P} achieves and $\epsilon$-approximate solution with $\mathcal{O}\left(\log(1/\epsilon)\right)$ samples. 
\end{proposition}

\begin{remark}
We stress that linear convergence with ZO2P was first established in \cite{mohammadi2020linear} for this problem and extended to continuous-time in~\cite{mohammadi2020learning, mohammadi2021convergence}. However, in Appendix \ref{appendix:proof_linear_ZO2P} we present an alternative and straightforward proof, one that relies simply on the upper bound of the expectation\footnote{Expectation  is taken with respect to  $U \stackrel{\text{i.i.d.}}{\sim} \mathcal{S}_r$ and $x_0  \stackrel{\text{i.i.d.}}{\sim} \mathcal{X}_0$.} of the estimated gradient, i.e.,  $\mc{E}\|\widetilde{\nabla}(K)\|^2_F$ (Lemma \ref{lemma:two_point_estimated_grad}) and does not involve proving that $\langle \widetilde{\nabla}C(K), \nabla C(K) \rangle \geq \mu_1 \|\nabla C(K)\|^2_F$, and $\|\widetilde{\nabla} C(K)\|^2_F \leq \mu_2 \|\nabla C(K)\|^2_F$ are satisfied with high probability, for $\mu_1, \mu_2 \in \mathbb{R}_{+}$  \cite[Section V]{mohammadi2020linear}.
\end{remark}


\begin{assumption}\label{assumption:Lipschitz_approx_gradient}  Let $\overline{g}(K) = \frac{d}{r^2}C(K+U)U$ be a single sample ZO1P gradient estimation with $U \stackrel{\text{i.i.d.}}{\sim} \mathcal{S}_r$. Then, for any two stabilizing policies $K$, $K^{\prime} \in \mathcal{G},$ we assume that
$$
    \mc{E}\|\overline{g}(K) - \overline{g}(K')\|_F \leq C_g\mc{E}\|K - K'\|_F.
$$
for some positive constant $C_g$.
\end{assumption} 
\begin{remark} Note that this assumption on the local smoothness of the estimated gradient is a standard requirement for variance-reduced algorithms, as established in \cite{khanduri2021stem, fang2018spider}. In the context of the LQR problem, this assumption has the same flavor as the local Lipschitz condition on the empirical cost function in \cite[Section 2]{malik2019derivative}.   
\end{remark}

Next, we present two auxiliary results that are instrumental in proving our main results. First, we control the bias in the zeroth-order gradient estimation (Lemma \ref{lemma:bias}) and establish a uniform bound for ZO2P estimated gradient (Lemma \ref{lemma:two_point_estimated_grad}).


\begin{lemma}  \label{lemma:bias} (Controlling the bias) Let $\widehat{\nabla}C(K)$ be the ZO1P or ZO2P gradient estimations evaluated at the stabilizing policy $K \in \mathcal{G}$. Then,
$$
    \mc{E}\|\nabla C(K) - \mc{E}\widehat{\nabla}C(K)\|^2_F \leq  \mathcal{B}(r) := \left(\barhgrad r\right)^2.
$$

\begin{proof} See Appendix \ref{appendix:proof_bias}. 
\end{proof}   
\end{lemma}


\begin{lemma} \label{lemma:two_point_estimated_grad}
Let $\widetilde{\nabla}(K)$ be the ZO2P gradient estimation. For any stabilizing policy $K \in \mathcal{G}$, it holds that 
$$
    \mc{E}\|\widetilde{\nabla}(K)\|^2_F \leq 8d^2\mathcal{B}(r)+ 2d^2\mc{E}\|\nabla C(K) \|^2_F.
$$

\begin{proof} 
See Appendix \ref{appendix:proof_lemma_two_point_estimated_grad}.
\end{proof}   
\end{lemma}

\subsection{Stability Analysis}

We now introduce the conditions on the number of samples $\{n_1,n_2\}$, step-size $\eta$ and smoothing radius $r$ to ensure that Algorithm \ref{algorithm:LQR_SVRPG} produces a stabilizing policy $K_{t+1}^{n+1}$ at each epoch $n \in \{0,\ldots,N-1\}$ and each $t \in \{0,\ldots,T-1\}$. 

\begin{theorem} (Per-iteration Stability) \label{theorem:stability}
Given  $K_0 \in \mathcal{G}$, suppose we set the number of outer and inner-loop samples such that satisfies $\{n_1, n_2\} \gtrsim \bar{h}_s\left(\frac{\psi}{6}, \delta\right)$, the step-size $\eta \lesssim  \frac{r^2\Delta_0}{\barhgrad d^2}$, and the smoothing radius
$$
r\leq \underline{h}_r\left(\frac{\psi}{6}\right):=\min \left\{\underline{h}_{\Delta}, \frac{1}{\bar{h}_{\text {cost }}}, \frac{\psi}{6 \bar{h}_{\text {grad }}}\right\},   
$$
with $\delta \in (0,1)$, $\psi := \sqrt{\frac{\lambda \Delta_0}{4}}$. Then, with probability $1-\delta$, it holds that $K^{n+1}_{t+1} \in \mathcal{G}$, for all $n$ and $t$.

\begin{proof}A detailed proof with the explicitly expression of $\bar{h}_s\left(\frac{\psi}{6}, \delta\right)$ is provided in Appendix \ref{appendix:proof_stability}. \\   
\end{proof}   
\end{theorem}

\noindent \textbf{Discussion:} We emphasize that, unlike the RL setting in \cite{papini2018stochastic, xu2020improved}, in the LQR optimal control problem, it is imperative to ensure the closed-loop stability of \eqref{eq:LTI_system} under $K^{n+1}_{t+1}$ for all $n \in \{0,\ldots,N-1\}$ and $t \in \{0,\ldots,T-1\}$. However, despite its dual-loop structure, demonstrating that $K^{n+1}_{t+1} \in \mathcal{G}$ throughout the iterations of Algorithm \ref{algorithm:LQR_SVRPG} can be achieved by following a similar approach as outlined in previous works without variance reduction~\cite{fazel2018global, gravell2020learning, malik2019derivative, mohammadi2020linear}. 

To this end, we first set the first iteration as the base case and demonstrate that as long as $K_0 \in \mathcal{G}$ (Assumption \ref{assumption:stabilizing_K0}), then $C(K^1_1) - C(K^*) \leq C(K_0) - C(K^*)$ holds true, indicating that $K^1_1 \in \mathcal{G}$. To establish this, we use the Lipschitz property of the cost function (Lemma \ref{lemma:Lipschitz_true_cost}), along with the gradient domination condition (Lemma \ref{lemma:PL_condition}), and the matrix Bernstein inequality \cite[Section 6]{tropp2012user}. The latter provides the necessary conditions on ${n_1, n_2}$ and $r$ to upper bound $\|\nabla C(K) - v^1_0\|_F \leq \psi$. The stability analysis is then completed by applying an induction step to this base case.

\subsection{Convergence Analysis}

We now proceed with our analysis to provide the necessary conditions on the number of samples $\{n_1, n_2\}$, smoothing radius $r$, step-size $\eta$, and total number of iterations $NT$ to ensure the global convergence of Algorithm \ref{algorithm:LQR_SVRPG}.

\begin{theorem} (Convergence Analysis) \label{theorem:convergence}
Suppose we select $n_2 \geq \max\left\{96d^2, \frac{\left(3C^2_g + 12\barhgrad^2d^2\right) T^2}{\barhgrad^2}\right\}$, and $\eta \leq  \frac{1}{4\barhgrad}$. Then, the policy $K_{out}$ returned by Algorithm \ref{algorithm:LQR_SVRPG} after $NT$ iterations enjoys the following property:
$$
    \mc{E}\left(C(K_{\text{out}}) - C(K^*)\right) \leq \Delta_0 \times  \left(1 - \frac{\eta \lambda}{16}\right)^{NT}  + \frac{\mathcal{B}(r)\phi}{\lambda n_2}.
$$
with $\phi = 120+ 192d^2$.

\begin{proof} 
Below we provide the proof strategy for this theorem. A detailed proof is presented in Appendix \ref{appendix:proof_convergence}. \\
\end{proof}   

\noindent \textbf{Proof Sketch:} Theorem \ref{theorem:convergence} is proved as follows:

\noindent 1) With the fact that $K^{n+1}_{t+1} \in \mathcal{G}$ for all $n \in \{0,\ldots,N-1\}$ and $t \in \{0,\ldots,T-1\}$ (Theorem \ref{theorem:stability}), along with Lemma \ref{lemma:Lipschitz_true_cost} and Young's inequality we can  write
    \begin{align}\label{eq:Lips_expectation}
        &\mc{E}\left(C(K^{n+1}_{t+1}) - C(K^{n+1}_{t})\right) \leq \frac{3\eta}{4} \mc{E}\|\nabla C(K^{n+1}_t) - v^{n+1}_t\|^2_F\notag \\
        &- \frac{\eta}{8}\mc{E}\|\nabla C(K^{n+1}_t)\|^2_F- \frac{\barhgrad}{2}\mc{E}\|K^{n+1}_{t+1} - K^{n+1}_t\|^2_F,
     \end{align}

\noindent 2) We control $\mc{E}\|\nabla C(K^{n+1}_t) - v^{n+1}_t\|^2_F$ in the above expression by decomposing it into bias and variance terms. In particular, we have:  biases from the inner and outer-loop estimations + variance of the ZO2P outer-loop estimation + ZO1P gradient estimation difference at $K^{n+1}_{t}$ and $\tilde{K}^n$. Both ZO1P and ZO2P biases are controlled in Lemma \ref{lemma:bias}. For the variance of the ZO2P gradient estimation we use Lemma \ref{lemma:two_point_estimated_grad} and for the ZO1P gradient difference term we assume local smoothness (Assumption \ref{assumption:Lipschitz_approx_gradient}). Thus, with $n_2 \geq 96d^2$, we have
     \begin{align*}
        &\mc{E}\|\nabla C(K^{n+1}_t) - v^{n+1}_t\|^2_F \leq \frac{\phi \eta \mathcal{B}(r)}{16n_2} + \tilde{\phi}\mc{E}\| K^{n+1}_t-\Tilde{K}^n\|^2_F\\
        &+\frac{1}{16}\mc{E}\|{\nabla}C(K_{t}^{n+1})\|^2_F, \text{ with } \tilde{\phi} = \frac{4}{3n_2}\left(\frac{3C^2_g}{2} + 6\barhgrad^2d^2\right).
     \end{align*}

\noindent 3) The proof is completed by using the PL condition (Lemma \ref{lemma:PL_condition}) and telescoping \eqref{eq:Lips_expectation} over outer and inner-loop iterations, with $n_2 \geq \frac{\left(3C^2_g + 12\barhgrad^2d^2\right) T^2}{\barhgrad^2}$, and $\eta \leq  \frac{1}{4\barhgrad}$.

\end{theorem}

\begin{corollary}(Sample Complexity)\label{corollary:sample_complexity}
Under the conditions of Theorem \ref{theorem:convergence}, and suppose we select the total number of iterations and smoothing radius according to 
$$
    NT \geq \frac{16 \log\left( 2\Delta_0/\epsilon\right)}{\eta \lambda},\;\ r \leq \sqrt{\frac{n_2 \lambda \epsilon}{2\phi\bar{h}^2_{\text{grad}}}},
$$
then Algorithm \ref{algorithm:LQR_SVRPG} achieves  $\mc{E}\left(C(K_{\text{out}}) - C(K^*)\right) \leq \epsilon$ with $\mathcal{O}\left(\log\left(1/\epsilon\right)^{3-2\beta}\right)$ cost queries. 

\begin{proof} The total number of cost queries required in Algorithm \ref{algorithm:LQR_SVRPG} is given by $\mathbb{S}_c :=  NTn_2 + Nn_1$. Therefore, since $n_1 = \mathcal{O}(1)$, the sample complexity of Algorithm \ref{algorithm:LQR_SVRPG} is dominated by the order of $NTn_2$. As a result, by setting $N=\mathcal{O}\left(\log\left(1/\epsilon\right)\right)^{\beta}$ and $T=\mathcal{O}\left(\log\left(1/\epsilon\right)\right)^{1-\beta}$, Algorithm \ref{algorithm:LQR_SVRPG} returns an $\epsilon$-approximate solution with $\mathcal{O}\left(\log\left(1/\epsilon\right)^{3-2\beta}\right)$ total number of cost queries.
\end{proof}   
\end{corollary}

\begin{corollary} (Oracle Complexity Reduction)\label{corollary:oracle_reduction}
Under the conditions of Theorem \ref{theorem:convergence} and Corollary \ref{corollary:sample_complexity}, it holds that Algorithm \ref{algorithm:LQR_SVRPG} achieves an $\epsilon$-approximate solution with a reduction of $\mathcal{O}\left(\log\left(1/\epsilon\right)\right)^{1-\beta}$ in the two-point cost queries when compared to Algorithm \ref{algorithm:LQR_ZO2P}, where $\beta \in (0,1)$.\\
\end{corollary}

\noindent \textbf{Discussion:} Similar to Corollary \ref{corollary:sample_complexity}, we set $N=\mathcal{O}\left(\log\left(1/\epsilon\right)\right)^{\beta}$ and $T=\mathcal{O}\left(\log\left(1/\epsilon\right)\right)^{1-\beta}$. Then, we observe that Algorithm \ref{algorithm:LQR_SVRPG}, with number of outer-loop samples $n_1 = \mathcal{O}(1)$, demands only $\mathcal{O}\left(\log\left(1/\epsilon\right)\right)^{\beta}$ two-point queries (i.e., the more resource-intensive cost queries to obtain) to achieve an $\epsilon$-approximate solution. This improves upon the two-point oracle complexity of Algorithm \ref{algorithm:LQR_ZO2P} by a factor of $\mathcal{O}\left(\log\left(1/\epsilon\right)\right)^{1-\beta}$. To verify this we simply note that our algorithm necessitates $\mathcal{N}_{\text{ZO2P}}= Nn_1= \mathcal{O}\left(\log\left(1/\epsilon\right)\right)^{\beta}$, whereas Algorithm \ref{algorithm:LQR_ZO2P} requires $\mathcal{N}_{\text{ZO2P}}=\mathcal{O}\left(\log\left(1/\epsilon\right)\right)$ two-point cost queries to attain $\mc{E}\left(C(K_{\text{out}}) - C(K^*)\right) \leq \epsilon$.

\section{Numerical Experiments} \label{sec:numerical_results}

Numerical experiments \footnote{Code for exact reproduction of the proposed experiments can be downloaded from \url{https://github.com/jd-anderson/LQR_SVRPG}} are now conducted to illustrate and evaluate the effectiveness of Algorithm \ref{algorithm:LQR_SVRPG}.  To ensure a fair comparison on the performance of the algorithms we set $x^\top_0 = [1 , 1 , 1]$ for computing the normalized cost gap between the current and optimal cost, namely,  $\frac{C(K_{l})-C(K^*)}{C(K_0)-C(K^*)}$, and $\mathcal{X}_0 \stackrel{d}{=}\mathcal{N}(0,I_{n_x})$ for the cost oracle generation.

Consider a unstable system with $n_x=3$ states and $n_u=1$ input, where the system and cost matrices are detailed in Appendix \ref{appendix:exp_additional_details}.  We set the initialization parameters of Algorithms \ref{algorithm:LQR_ZO2P} and \ref{algorithm:LQR_SVRPG} as follows: 1) $r= 1\times 10^{-4}$, $n_1 = 50$, $\eta = 1\times 10^{-4}$. 2) 
$r_{\text{in}}= 5\times 10^{-2}$, $r_{\text{out}}=1\times 10^{-4}$, $n_1 = 50$, $n_2 = 25$, $N = 125$, $T=4$,  $\eta = 1\times 10^{-4}$.

\begin{figure}
     \centering
         \includegraphics[width=0.5\textwidth]{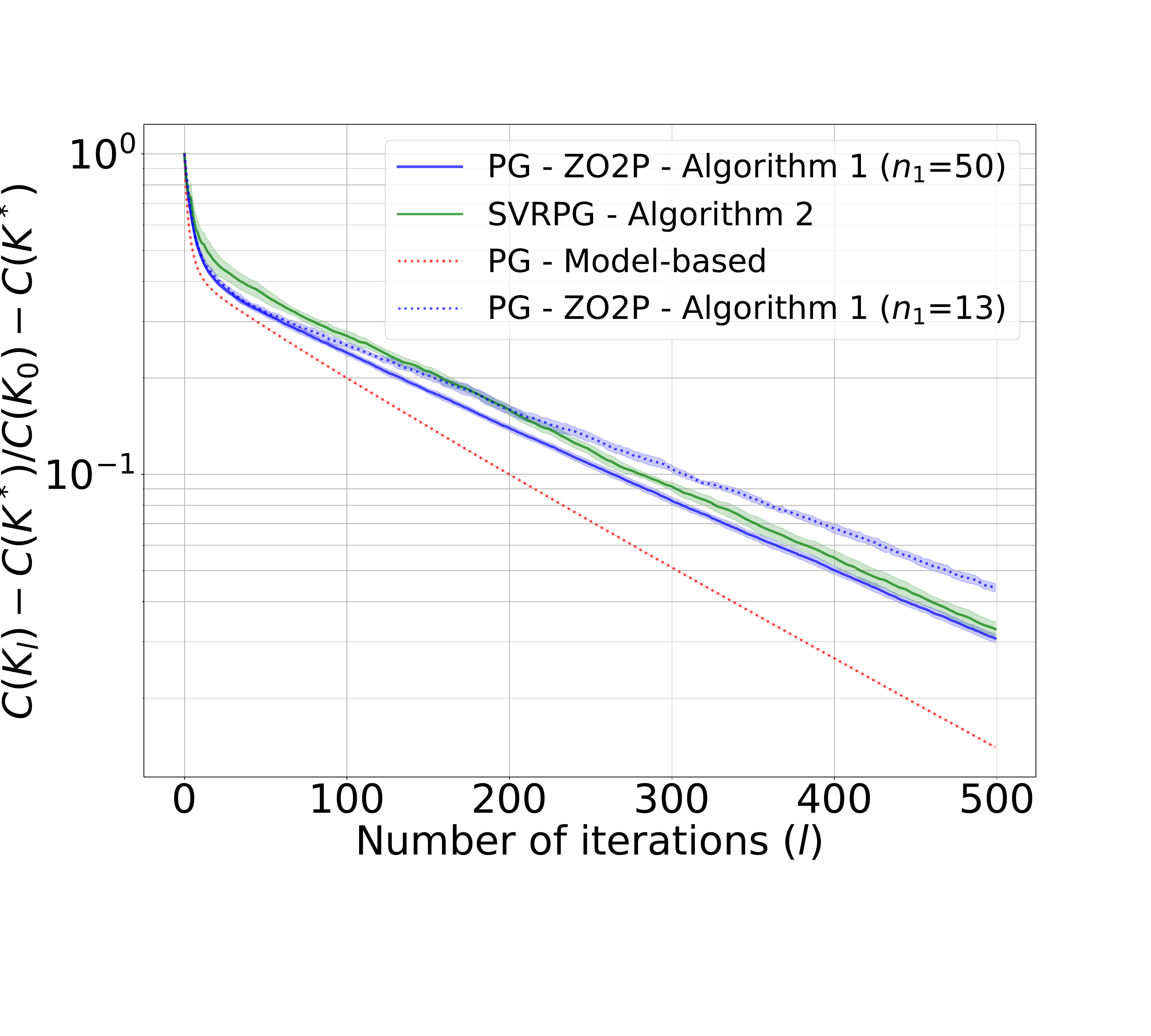}
        \vspace{-1.1cm}\caption{ Normalized gap between the current and optimal cost with respect to the iteration count.} 
        \label{fig:numerics}
\end{figure}

Figure \ref{fig:numerics} demonstrates the convergence of Algorithms \ref{algorithm:LQR_ZO2P} and \ref{algorithm:LQR_SVRPG}. It also includes the result for the policy gradient descent under the model-based setting. The latter highlights the limit of how well the PG algorithms discussing in this work can do without knowing the system model. 

The figure shows that both Algorithms \ref{algorithm:LQR_ZO2P} and \ref{algorithm:LQR_SVRPG} achieve an equivalent convergence performance for the specified parameters. We emphasize that Algorithms \ref{algorithm:LQR_ZO2P} and \ref{algorithm:LQR_SVRPG} use $\mathcal{S}_c = 50000$ and $\mathcal{S}_c = 37500$  cost queries, respectively, to attain $\epsilon = 3\times 10^{-2}$. Moreover, in terms of two-point queries, Algorithm \ref{algorithm:LQR_SVRPG} necessitates only $\mathcal{N}_{\text{ZO2P}}= Nn_1= 6250$, whereas Algorithm \ref{algorithm:LQR_ZO2P} is entirely reliant on two-point queries, requiring $25000$ to achieve the same accuracy as shown in the figure. The figure also shows that the performance of Algorithm \ref{algorithm:LQR_ZO2P} degrades when the number of two-point queries decreases to $6500$. This demonstrates that with our SVRPG approach we are able to effectively reduce the two-point oracle complexity for solving the model-free LQR problem.

\section{Conclusions and Future Work} \label{sec:conclusions}

We proposed an oracle efficient algorithm to solve the model-free LQR problem. Our approach combines a SVRPG-based approach with a mixed zeroth-order gradient estimation scheme. This mixed gradient estimation yields a reduction in the number of two-point cost queries required to achieve an $\epsilon$-approximate solution since the more resource-expensive queries are now required less frequently. We proved that our approach improves by a factor of $\mathcal{O}\left(\log\left(1/\epsilon\right)\right)^{1-\beta}$ two-point query information upon the standard ZO2P gradient estimation method. Future work will involve exploring loopless variants and recursive momentum-based approaches to further reduce the two-point oracle complexity required to solve the model-free LQR problem.

\bibliography{bibliography.bib}

\begin{thebibliography}{10}
\providecommand{\url}[1]{#1}
\csname url@samestyle\endcsname
\providecommand{\newblock}{\relax}
\providecommand{\bibinfo}[2]{#2}
\providecommand{\BIBentrySTDinterwordspacing}{\spaceskip=0pt\relax}
\providecommand{\BIBentryALTinterwordstretchfactor}{4}
\providecommand{\BIBentryALTinterwordspacing}{\spaceskip=\fontdimen2\font plus
\BIBentryALTinterwordstretchfactor\fontdimen3\font minus \fontdimen4\font\relax}
\providecommand{\BIBforeignlanguage}[2]{{%
\expandafter\ifx\csname l@#1\endcsname\relax
\typeout{** WARNING: IEEEtran.bst: No hyphenation pattern has been}%
\typeout{** loaded for the language `#1'. Using the pattern for}%
\typeout{** the default language instead.}%
\else
\language=\csname l@#1\endcsname
\fi
#2}}
\providecommand{\BIBdecl}{\relax}
\BIBdecl

\bibitem{ziemann2022policy}
I.~Ziemann, A.~Tsiamis, H.~Sandberg, and N.~Matni, ``How are policy gradient methods affected by the limits of control?'' in \emph{2022 IEEE 61st Conference on Decision and Control (CDC)}.\hskip 1em plus 0.5em minus 0.4em\relax IEEE, 2022, pp. 5992--5999.

\bibitem{fazel2018global}
M.~Fazel, R.~Ge, S.~Kakade, and M.~Mesbahi, ``Global convergence of policy gradient methods for the linear quadratic regulator,'' in \emph{International conference on machine learning}.\hskip 1em plus 0.5em minus 0.4em\relax PMLR, 2018, pp. 1467--1476.

\bibitem{gravell2020learning}
B.~Gravell, P.~M. Esfahani, and T.~Summers, ``Learning optimal controllers for linear systems with multiplicative noise via policy gradient,'' \emph{IEEE Transactions on Automatic Control}, vol.~66, no.~11, pp. 5283--5298, 2020.

\bibitem{malik2019derivative}
D.~Malik, A.~Pananjady, K.~Bhatia, K.~Khamaru, P.~Bartlett, and M.~Wainwright, ``Derivative-free methods for policy optimization: Guarantees for linear quadratic systems,'' in \emph{The 22nd international conference on artificial intelligence and statistics}.\hskip 1em plus 0.5em minus 0.4em\relax PMLR, 2019, pp. 2916--2925.

\bibitem{mohammadi2020linear}
H.~Mohammadi, M.~Soltanolkotabi, and M.~R. Jovanovi{\'c}, ``On the linear convergence of random search for discrete-time lqr,'' \emph{IEEE Control Systems Letters}, vol.~5, no.~3, pp. 989--994, 2020.

\bibitem{hu2023toward}
B.~Hu, K.~Zhang, N.~Li, M.~Mesbahi, M.~Fazel, and T.~Ba{\c{s}}ar, ``Toward a theoretical foundation of policy optimization for learning control policies,'' \emph{Annual Review of Control, Robotics, and Autonomous Systems}, vol.~6, pp. 123--158, 2023.

\bibitem{perdomo2021stabilizing}
J.~Perdomo, J.~Umenberger, and M.~Simchowitz, ``Stabilizing dynamical systems via policy gradient methods,'' \emph{Advances in neural information processing systems}, vol.~34, pp. 29\,274--29\,286, 2021.

\bibitem{spall2005introduction}
J.~C. Spall, \emph{Introduction to stochastic search and optimization: estimation, simulation, and control}.\hskip 1em plus 0.5em minus 0.4em\relax John Wiley \& Sons, 2005.

\bibitem{peters2008reinforcement}
J.~Peters and S.~Schaal, ``Reinforcement learning of motor skills with policy gradients,'' \emph{Neural networks}, vol.~21, no.~4, pp. 682--697, 2008.

\bibitem{wang202model}
H.~Wang, L.~F. Toso, A.~Mitra, and J.~Anderson, ``{Model-free Learning with Heterogeneous Dynamical Systems: A Federated LQR Approach},'' \emph{arXiv preprint arXiv:2308.11743}, 2023.

\bibitem{tang2023zeroth}
Y.~Tang, Z.~Ren, and N.~Li, ``Zeroth-order feedback optimization for cooperative multi-agent systems,'' \emph{Automatica}, vol. 148, p. 110741, 2023.

\bibitem{johnson2013accelerating}
R.~Johnson and T.~Zhang, ``Accelerating stochastic gradient descent using predictive variance reduction,'' \emph{Advances in neural information processing systems}, vol.~26, 2013.

\bibitem{papini2018stochastic}
M.~Papini, D.~Binaghi, G.~Canonaco, M.~Pirotta, and M.~Restelli, ``Stochastic variance-reduced policy gradient,'' in \emph{International conference on machine learning}.\hskip 1em plus 0.5em minus 0.4em\relax PMLR, 2018, pp. 4026--4035.

\bibitem{mohammadi2019global}
H.~Mohammadi, A.~Zare, M.~Soltanolkotabi, and M.~R. Jovanovi{\'c}, ``Global exponential convergence of gradient methods over the nonconvex landscape of the linear quadratic regulator,'' in \emph{2019 IEEE 58th Conference on Decision and Control (CDC)}.\hskip 1em plus 0.5em minus 0.4em\relax IEEE, 2019, pp. 7474--7479.

\bibitem{mohammadi2020random}
H.~Mohammadi, M.~Soltanolkotabi, and M.~R. Jovanovic, ``Random search for learning the linear quadratic regulator,'' in \emph{2020 American Control Conference (ACC)}.\hskip 1em plus 0.5em minus 0.4em\relax IEEE, 2020, pp. 4798--4803.

\bibitem{mohammadi2021convergence}
H.~Mohammadi, A.~Zare, M.~Soltanolkotabi, and M.~R. Jovanovi{\'c}, ``Convergence and sample complexity of gradient methods for the model-free linear--quadratic regulator problem,'' \emph{IEEE Transactions on Automatic Control}, vol.~67, no.~5, pp. 2435--2450, 2021.

\bibitem{roux2012stochastic}
N.~Roux, M.~Schmidt, and F.~Bach, ``A stochastic gradient method with an exponential convergence \_rate for finite training sets,'' \emph{Advances in neural information processing systems}, vol.~25, 2012.

\bibitem{defazio2014saga}
A.~Defazio, F.~Bach, and S.~Lacoste-Julien, ``{SAGA}: A fast incremental gradient method with support for non-strongly convex composite objectives,'' \emph{Advances in neural information processing systems}, vol.~27, 2014.

\bibitem{xu2020improved}
P.~Xu, F.~Gao, and Q.~Gu, ``An improved convergence analysis of stochastic variance-reduced policy gradient,'' in \emph{Uncertainty in Artificial Intelligence}.\hskip 1em plus 0.5em minus 0.4em\relax PMLR, 2020, pp. 541--551.

\bibitem{liu2020improved}
Y.~Liu, K.~Zhang, T.~Basar, and W.~Yin, ``An improved analysis of (variance-reduced) policy gradient and natural policy gradient methods,'' \emph{Advances in Neural Information Processing Systems}, vol.~33, pp. 7624--7636, 2020.

\bibitem{ji2019improved}
K.~Ji, Z.~Wang, Y.~Zhou, and Y.~Liang, ``Improved zeroth-order variance reduced algorithms and analysis for nonconvex optimization,'' in \emph{International conference on machine learning}.\hskip 1em plus 0.5em minus 0.4em\relax PMLR, 2019, pp. 3100--3109.

\bibitem{liu2018zeroth}
S.~Liu, B.~Kailkhura, P.-Y. Chen, P.~Ting, S.~Chang, and L.~Amini, ``Zeroth-order stochastic variance reduction for nonconvex optimization,'' \emph{Advances in Neural Information Processing Systems}, vol.~31, 2018.

\bibitem{hewer1971iterative}
G.~Hewer, ``An iterative technique for the computation of the steady state gains for the discrete optimal regulator,'' \emph{IEEE Transactions on Automatic Control}, vol.~16, no.~4, pp. 382--384, 1971.

\bibitem{mohammadi2020learning}
H.~Mohammadi, M.~R. Jovanovic, and M.~Soltanolkotabi, ``Learning the model-free linear quadratic regulator via random search,'' in \emph{Learning for Dynamics and Control}.\hskip 1em plus 0.5em minus 0.4em\relax PMLR, 2020, pp. 531--539.

\bibitem{khanduri2021stem}
P.~Khanduri, P.~Sharma, H.~Yang, M.~Hong, J.~Liu, K.~Rajawat, and P.~Varshney, ``Stem: A stochastic two-sided momentum algorithm achieving near-optimal sample and communication complexities for federated learning,'' \emph{Advances in Neural Information Processing Systems}, vol.~34, pp. 6050--6061, 2021.

\bibitem{fang2018spider}
C.~Fang, C.~J. Li, Z.~Lin, and T.~Zhang, ``Spider: Near-optimal non-convex optimization via stochastic path-integrated differential estimator,'' \emph{Advances in neural information processing systems}, vol.~31, 2018.

\bibitem{tropp2012user}
J.~A. Tropp, ``User-friendly tail bounds for sums of random matrices,'' \emph{Foundations of computational mathematics}, vol.~12, pp. 389--434, 2012.

\bibitem{gao2018information}
X.~Gao, B.~Jiang, and S.~Zhang, ``On the information-adaptive variants of the admm: an iteration complexity perspective,'' \emph{Journal of Scientific Computing}, vol.~76, pp. 327--363, 2018.

\end{thebibliography}
\bibliographystyle{IEEEtran}

\onecolumn
\appendix

\allowdisplaybreaks

\subsection{Auxiliary Definitions}\label{appendix:auxiliary_def}

Before delving into the proofs, let us begin by defining the following quantities:
\begin{align*}
&\bar{h}_{\text{grad}} := \sup_{K \in \mathcal{G}} {h}_{\text{grad}}(K), \ \bar{h}_{\text{cost}} := \sup_{K \in \mathcal{G}} {h}_{\text{cost}}(K),\ 
 \underline{h}_{\Delta} := \inf_{K \in \mathcal{G}} h_{\Delta}(K),\notag\\
& C_{\max} := \max_{K \in \mathcal{G}} C(K),\;\ \|Q_K\|_{\max} := \max_{K \in \mathcal{G}} \|Q + K^\top RK\|_F.
\end{align*}

\subsection{Proof of Proposition \ref{prop:linear_ZO2P}} \label{appendix:proof_linear_ZO2P}

We start our proof with the local smoothness of the cost function (Lemma \ref{lemma:Lipschitz_true_cost}) to write

\begin{align*}
   \mc{E}\left(C(K_{l+1}) - C(K_l)\right) &\leq -\eta \mc{E}\langle \nabla C(K_l), \widetilde{\nabla} C(K_l) \rangle + \frac{\bar{h}_{\text{grad}}\eta^2}{2}\mc{E}\| \widetilde{\nabla}C(K_l)\|^2_F\\
   & \leq -\eta \mc{E}\| \nabla C(K_l) \|^2_F + \frac{\eta}{2}\mc{E}\| \nabla C(K_l) \|^2_F+ \frac{\eta}{2}\|\mc{E}\widetilde{\nabla}C(K_l) - \nabla C(K_l)\|^2_F  + \frac{\bar{h}_{\text{grad}}\eta^2}{2}\mc{E}\| \widetilde{\nabla}C(K_l)\|^2_F\\
   &\stackrel{(i)}{\leq}  -\frac{\eta}{2}\mc{E}\| \nabla C(K_l) \|^2_F + \frac{\eta \left(\barhgrad r\right)^2}{2} + \frac{\bar{h}_{\text{grad}}\eta^2}{2}\mc{E}\| \widetilde{\nabla}C(K_l)\|^2_F\\
   &\stackrel{(ii)}{\leq}  -\frac{\eta}{4}\mc{E}\| \nabla C(K_l) \|^2_F + \frac{3\eta\left(\barhgrad r\right)^2}{2}
\end{align*}
where $(i)$ is the due to the upper bound on the estimation bias (Lemma \ref{lemma:bias}), and $(ii)$ follows from the bound on $\mc{E}\|\widetilde{\nabla}C(K_l)\|^2_F$ in Lemma \ref{lemma:two_point_estimated_grad} and the selection of the step-size $\eta \leq \frac{1}{4\bar{h}_{\text{grad}}d^2}$. Let us now exploit the PL condition (Lemma \ref{lemma:PL_condition}) to obtain
\begin{align*}
   &\mc{E}\left(C(K_{l+1}) - C(K^*)\right) \leq \tilde{\gamma} \mc{E}\left(C(K_l) - C(K^*)\right)+ \frac{3\eta\left(\barhgrad r\right)^2}{2},
\end{align*}
with $\tilde{\gamma}=1 -\frac{\eta\lambda}{4}$ which can be telescoped over $l \in \{0,1,\ldots,L-1\}$ to obtain 
\begin{align*}
   &\mc{E}\left(C(K^{L}) - C(K^*)\right) \leq \tilde{\gamma}^L\Delta_0+ \frac{6(\barhgrad r)^2}{\lambda},
\end{align*}
with $\Delta_0 = C(K_0) - C(K^*)$. Thus, by selecting the total number of iterations $L$ and smoothing radius $r$ according to
\begin{align*}
   &L \geq \frac{4}{\lambda \eta}\log\left(\frac{2\Delta_0}{\epsilon}\right) = \mathcal{O}\left(\log\left(1/\epsilon\right)\right), \;\ r \leq \sqrt{\frac{\epsilon \lambda }{12\bar{h}^2_{\text{grad}}}} = \mathcal{O}\left( \sqrt{\epsilon}\right)
\end{align*}\\
with number of samples $n_1 = \mathcal{O}(1)$, we thus obtain $\mc{E}\left(C(K^{L}) - C(K^*)\right) \leq \epsilon$, which completes the proof.

\subsection{Proof of Lemma \ref{lemma:two_point_estimated_grad}} \label{appendix:proof_lemma_two_point_estimated_grad}

By definition we have

\allowdisplaybreaks
\begin{align*}
    \mc{E}\| \widetilde{\nabla}C(K)\|^2_F &=  \frac{d^2}{4r^4}\mc{E}\left\|\frac{1}{n_1}\sum_{i=1}^{n_1}  \left( C(K + U_i) - C(K - U_i) \right)U_i \right\|^2_F\\
    & \stackrel{(i)}{\leq} \frac{d^2n_1}{4n^2_1r^2}\sum_{i=1}^{n_1}\mc{E}\left( C(K + U_i) - C(K - U_i) \right)^2\\
    & = \frac{d^2}{4n_1r^2}\sum_{i=1}^{n_1}\mc{E}\left( C(K + U_i) - C(K- U_i)-\langle \nabla C(K), 2U_i \rangle  + \langle \nabla C(K), 2U_i \rangle \right)^2\\
    &\stackrel{(ii)}{\leq} \frac{d^2}{4n_1r^2}\sum_{i=1}^{n_1} 2\mc{E}\left( C(K + U_i) - C(K - U_i)-\langle \nabla C(K), 2U_i \rangle \right)^2  + 2\mc{E}\left(\langle \nabla C(K), 2U_i \rangle \right)^2\\
    &\leq \frac{d^2}{4n_1r^2}\sum_{i=1}^{n_1} 2\mc{E}\left( C(K + U_i) - C(K - U_i)- \langle \nabla C(K), 2U_i \rangle \right)^2+ 8r^2\mc{E}\|\nabla C(K) \|^2_F\\
    &= \frac{d^2}{4n_1r^2}\sum_{i=1}^{n_1}\left( 2\mc{E}\left( C(K + U_i) - C(K - U_i) - \langle \nabla C(K), 2U_i \rangle -\langle \nabla C(K - U_i), 2U_i \rangle + \langle \nabla C(K - U_i), 2U_i \rangle   \right)^2\right.\\
   & \left.+ 8r^2\mc{E}\|\nabla C(K) \|^2_F\right)\\
    &\leq \frac{d^2}{4n_1r^2}\sum_{i=1}^{n_1}\left( 4\mc{E}\left( C(K + U_i) - C(K - U_i)- \langle \nabla C(K - U_i), 2U_i \rangle\right)^2+4\mc{E}\left(\langle \nabla C(K - U_i) - \nabla C(K), 2U_i \rangle   \right)^2\right.\\
    &\left.+ 8r^2\mc{E}\|\nabla C(K) \|^2_F\right)\\
    &\stackrel{(iii)}{\leq}  \frac{d^2}{4n_1r^2}\sum_{i=1}^{n_1} 16\left(\barhgrad r\right)^2 +4\mc{E}\left(\langle \nabla C(K - U_i)- \nabla C(K), 2U_i \rangle   \right)^2+ 8r^2\mc{E}\|\nabla C(K) \|^2_F\\
    &\stackrel{(iv)}{\leq} \frac{d^2}{4n_1r^2}\sum_{i=1}^{n_1} 32\left(\barhgrad r\right)^2+ 8r^2\mc{E}\|\nabla C(K) \|^2_F\\
    &= 8d^2\left(\barhgrad r\right)^2+ 2d^2\mc{E}\|\nabla C(K) \|^2_F,
\end{align*}
where $(i)$ is due to $\|U_i\|_F= r$, $(ii)$ follows from Young's inequality, $(iii)$ and $(iv)$ are due to the fact that the cost function is $\barhgrad$-Lipschitz.

\begin{remark}
A similar result is presented in Proposition 7.6 of \cite{gao2018information} for the stochastic gradient ADMM approach.
\end{remark}

\subsection{Proof of Lemma \ref{lemma:bias}} \label{appendix:proof_bias}

The proof for this lemma is presented in \cite{malik2019derivative}. For the self-completeness of our work we revisit it below. To prove this lemma we note that

\begin{align*}
    \mc{E}\|\nabla C(K) - \mc{E}\widehat{\nabla}C(K)\|^2_F &=  \mc{E}\|\nabla C(K) - {\nabla}\mc{E}C(K + U)\|^2_F\\
    &=\mc{E}\|\mc{E}\left(\nabla C(K) - {\nabla}C(K + U)\right)\|^2_F\\
    &\stackrel{(i)}{\leq} \mc{E}\mc{E}\|\nabla C(K) - {\nabla}C(K + U)\|^2_F\\
   &\stackrel{(ii)}{\leq} \mathcal{B}(r):=(\bar{h}_{\text{grad}}r)^2,
\end{align*}
where $(i)$ is due to Jensen's inequality and $(ii)$ follows from the gradient's Lipschitz property. It is also worth observing that the above result is satisfied for both ZO1P and ZO2P since
\begin{align*}
   \mc{E}\left[\frac{d}{2r^2}\left(C(K+U) - C(K-U)\right)U\right] &=\mc{E}\left[\frac{d}{2r^2}C(K+U)U\right] + \mc{E}\left[\frac{d}{2r^2}C(K-U)(-U)\right]\\
   &\stackrel{(i)}{=} \mc{E}\left[\frac{d}{2r^2}C(K+U)U\right] + \mc{E}\left[\frac{d}{2r^2}C(K+U)U\right]\\
   &=\mc{E}\left[\frac{d}{r^2}C(K+U)U\right],
\end{align*}
where $(i)$ follows from the symmetry of the distribution of $U_i$ on the sphere $\mathcal{S}_r$ of radius $r$.

\subsection{Proof of Theorem \ref{theorem:stability} (\textbf{Stability Analysis})} \label{appendix:proof_stability}

Now we proceed with our analysis to prove that $K^{n+1}_{t+1} \in \mathcal{G}$ for any iteration $t \in \{0,1,\ldots,T-1\}$ and  $n \in \{0,1,\ldots,N-1\}$ of Algorithm \ref{algorithm:LQR_SVRPG}, given that $K_0 \in \mathcal{G}$. Without loss of generality we set the base case (first iteration $t=0$, $  n=0$) and use an induction step to prove that $K^{n+1}_{t+1}$ is stabilizing for any other $t \in [T-1]$, and  $n \in [N-1]$. For this purpose, we start by exploiting the local smoothness of the cost function to write
\begin{align*}
  C(K^{1}_{1}) - C(K_{0}) &\leq \langle \nabla C(K_0), K^1_1 - K_0 \rangle + \frac{\bar{h}_{\text{grad}}}{2}\|K^1_1 - K_0\|^2_F,
\end{align*}
since $K^{1}_0 = K_0$. Now, we can control $\|K^1_1 - K_0\|_F$ by first writing
\begin{align*}
    \|K^1_1 - K_0\|^2_F &= \|\eta v^1_0\|^2_F = \eta^2 \left\|\widetilde{\nabla}C(K_0)+ \overline{\nabla}C(K_0) - \overline{\nabla} C(K_0)\right\|^2_F\\
    &\leq  \eta^2\left(2\|\widetilde{\nabla}C(K_0) \|_{F}  + 4\| \overline{\nabla}C(K_0)\|_F + 4\|\overline{\nabla} C(K_0) \|_F\right)
\end{align*}
where the above is due to triangle inequality. In addition, we have
\begin{align*}
    \|\widetilde{\nabla}C(K_0)\|^2_{F} &\stackrel{(i)}{\leq}  \frac{d^2}{4r^4n_1} \sum_{i=1}^{n_1} \left({C}(K_0 + U_i) - {C}(K_0 - U_i) \right)^2\|U_i\|^2_F\\ 
    &\stackrel{(ii)}{\leq} \frac{d^2}{2r^2n_1} \sum_{i=1}^{n_1} \left({C}(K_0 + U_i) \right)^2 + \left({C}(K_0 - U_i) \right)^2,
\end{align*}
where $(i)$ is due to Jensen's inequality, and $(ii)$ follows from Young's inequality and $\|U_i\|_F= r$. Thus, by using the local Lipschitz property of the cost function, i.e., 
\begin{align*}
  &|C(K_0 + U_i) - C(K_0)|\leq \bar{h}_{\text{cost}}C(K_0)\|U_i\|_F \leq \bar{h}_{\text{cost}}C(K_0)r,\\
  &|C(K_0 - U_i) - C(K_0)|\leq \bar{h}_{\text{cost}}C(K_0)\|U_i\|_F \leq \bar{h}_{\text{cost}}C(K_0)r,
\end{align*}
we have $C(K_0 + U_i) \leq 2C(K_0)$ and $C(K_0 - U_i) \leq 2C(K_0)$  for any smoothing radius that satisfies $r \leq \frac{1}{\bar{h}_{\text{cost}}}$. Then, we obtain,
\begin{align*}
    \|\widetilde{\nabla}C(K_0) \|^2_{F} \leq  \frac{4d^2{C}^2(K_0)}{r^2}.
\end{align*}

Similarly, $\| \overline{\nabla}C(K_0)\|_F$ and $\|\overline{\nabla} C(K_0)\|_F $ are upper bounded by
\begin{align*}
    \|\overline{\nabla}C(K_0) \|^2_{F} &\leq  \frac{d^2}{r^2 n_2} \sum_{i=1}^{n_2} {C}^2(K_0 + U_i)\|U_i\|^2_F\leq \frac{4d^2{C}^2(K_0)}{r^2},
\end{align*}
and
\begin{align*}
    \|{\overline{\nabla}}C(K_0) \|^2_{F}  &\leq  \frac{d^2}{r^2n_2} \sum_{i=1}^{n_2} {C}^2(K_0 + U_i)\|U_i\|^2_F\leq \frac{4d^2{C}^2(K_0)}{r^2},
\end{align*}
which yields
\begin{align*}
    \|K^1_1 - K_0\|^2_F &\leq \frac{40d^2{C}^2(K_0)\eta^2}{r^2}= \frac{\eta^2 c_1}{r^2},
\end{align*}
where $c_1 = 40d^2{C}^2(K_0)$. Therefore, we obtain
\begin{align*}
  C(K^{1}_{1}) - C(K_{0}) &\leq \langle \nabla C(K_0), -\eta v^1_0 \rangle + \frac{\eta^2\bar{h}_{\text{grad}}c_1}{2r^2}\\
  &=-\eta\langle \nabla C(K_0), \nabla C(K_0)+v^1_0 -\nabla C(K_0)\rangle+ \frac{\eta^2\bar{h}_{\text{grad}}c_1}{2r^2}\\
  &= -\eta \langle \nabla C(K_0), \nabla C(K_0)\rangle + \eta \langle \nabla C(K_0),\nabla C(K_0) -v^1_0\rangle+ \frac{\eta^2 \bar{h}_{\text{grad}}c_1}{2r^2}\\
  & = - \eta \|\nabla C(K_0)\|^2_F + \eta \langle \nabla C(K_0),\nabla C(K_0) -v^1_0\rangle+ \frac{\eta^2\bar{h}_{\text{grad}}c_1}{2r^2}\\
  &\stackrel{(i)}{\leq}- \frac{\eta}{2}\|\nabla C(K_0)\|^2_F + \frac{\eta}{2}\|\nabla C(K_0) -v^1_0\|^2_F+ \frac{\eta^2\bar{h}_{\text{grad}}c_1}{2r^2},
\end{align*}
where $(i)$ is due to Young's inequality on the second term. Thus, by using the PL condition (Lemma \ref{lemma:PL_condition}) we obtain

\begin{align*}
  C(K^{1}_{1}) - C(K^*) &\leq \left(1 -\frac{\eta\lambda}{2}\right)\left(C(K_0) - C(K^*) \right)+ \frac{\eta}{2}\|\nabla C(K_0) -v^1_0\|^2_F
  + \frac{\eta^2\bar{h}_{\text{grad}}c_1}{2r^2}\\
  &\stackrel{(i)}{\leq} \left(1 -\frac{\eta\lambda}{4}\right)\left(C(K_0) - C(K^*) \right)+ \frac{\eta}{2}\|\nabla C(K_0) -v^1_0\|^2_F,
\end{align*}
where $(i)$ follows from the selection of the step-size according to $\eta \leq  \frac{ r^2 \lambda \left(C(K_0) - C(K^*) \right)}{80\barhgrad d^2 C^2(K_0)}$. Let us now proceed with our analysis by controlling $\|\nabla C(K_0) -v^1_0\|^2_F$ such that  
\begin{align*}
    \|\nabla C(K_0) -v^1_0\|_F \leq \psi:= \sqrt{\frac{\lambda}{4}\left(C(K_0) - C(K^*) \right)},
\end{align*}
then we can write
\begin{align*}
    &\|\nabla C(K_0) -v^1_0\|_F = \left\|\nabla C(K_0) - \widetilde{\nabla}C(K_0) - \overline{\nabla}C(K_0) + \overline{\nabla} C(K_0)\right\|_F\\
    &=\left\|\left(\nabla C(K_0) - \widetilde{\nabla}C(K_0)\right) +\left(\nabla C(K_0)-\overline{\nabla}C(K_0) \right)+ \left(\overline{\nabla} C(K_0) -\nabla C(K_0)  \right)\right\|_F\\
    &\leq \underbrace{\left\|\nabla C(K_0) - \widetilde{\nabla}C(K_0) \right\|_F}_{\leq \frac{\psi}{3}} +\underbrace{\left\|\nabla C(K_0)-\overline{\nabla}C(K_0)\right\|_F}_{\leq \frac{\psi}{3}}+ \underbrace{\left\|\overline{\nabla} C(K_0) -\nabla C(K_0) \right\|_F}_{\leq \frac{\psi}{3}},
\end{align*}
which can controlled by using the matrix Bernstein inequality as in \cite[Lemma B.6]{gravell2020learning}. Thus, with probability $1-\delta$, with $\delta \in (0,1)$, it holds that $ \|\nabla C(K_0) -v^1_0\|_F \leq \psi$ for any inner and outer-loop samples and smoothing radius that satisfies
\begin{align}\label{eq:smoothing_rad_rollouts}
r \leq \underline{h}_r\left(\frac{\psi}{6}\right),
\{n_1,n_2\} \geq \bar{h}_s\left(\frac{\psi}{6}, \delta\right),
\end{align}
where 
\begin{align*}
\underline{h}_r\left(\frac{\psi}{6}\right):=\min \left\{\underline{h}_{\Delta}, \frac{1}{\bar{h}_{\text {cost }}}, \frac{\psi}{6 \bar{h}_{\text {grad }}}\right\},   
\end{align*}
\begin{align*}
\bar{h}_s\left(\frac{\psi}{6}, \delta\right) :=\frac{72 d_{\min}}{\psi^2}\left(\sigma_r^2+\frac{B_r \psi}{18 \sqrt{d_{\min}}}\right) \log \left[\frac{n_x+n_u}{\delta}\right],
\end{align*}
with $d_{\min}=\min (n_x, n_u)$,
\begin{align*}
B_r  :=\frac{2 d C_{\text{max}}}{r }+\frac{\psi}{6}+\bar{h}_1, \;\ \sigma_r^2  :=\left(\frac{2 d  C_{\text{max}}}{r }\right)^2+\left(\frac{\psi}{6}+\bar{h}_1\right)^2,
\end{align*}
and 
\begin{align*}
    C_{\max} := \max_{K \in \mathcal{G}} C(K),\;\ \|Q_K\|_{\max} := \max_{K \in \mathcal{G}} \|Q + K^TRK\|.
\end{align*}
which implies
\begin{align*}
  C(K^{1}_{1}) - C(K^*) &\leq \left(1 - \frac{\eta\lambda}{8}\right)\left(C(K_0) - C(K^*) \right),
\end{align*}
thus for any step-size $\eta>0$ we obtain 
\begin{align*}
  C(K^{1}_{1}) - C(K^*) &\leq C(K_0) - C(K^*) \to C(K^{1}_{1}) \leq C(K_0),
\end{align*}
which implies  $K^{1}_{1} \in \mathcal{G}$. Hence, the proof is completed by applying a induction step for all iterations with $t\in [T-1]$ and $n\in[N-1]$ in Algorithm \ref{algorithm:LQR_SVRPG}. To then establish that $K^{n+1}_{t+1} \in \mathcal{G}$ when the smoothing radius $r$, the number of samples $n_1$, and $n_2$ satisfy \eqref{eq:smoothing_rad_rollouts}.

\subsection{Proof of Theorem \ref{theorem:convergence} (\textbf{Convergence Analysis})} \label{appendix:proof_convergence}

Let us start by recalling the updating rule of Algorithm \ref{algorithm:LQR_SVRPG} 

\begin{align*}
    K^{n+1}_{t+1} = K^{n+1}_t - \eta v^{n+1}_t,
\end{align*}
with
\begin{align}\label{eq:stochastic_gradient}
v^{n+1}_t=\widetilde{\nabla}C(\tilde{K}^n) + \overline{\nabla} C(K^{n+1}_t) - \overline{\nabla} C(\Tilde{K}^n),
\end{align}
where $\widetilde{\nabla}C(\tilde{K}^n)$ corresponds to the ZO2P gradient estimation in the outer-loop and $\overline{\nabla}C(K^{n+1}_t)$, ${\overline{\nabla}}C(\tilde{K}^n)$ are the ZO1P estimation in the inner-loop. Now, we use the local smoothness property of the cost function (Lemma \ref{lemma:Lipschitz_true_cost}) to write
\begin{align*}
C(K^{n+1}_{t+1}) - C(K^{n+1}_{t}) &\leq \langle \nabla C(K^{n+1}_t), K^{n+1}_{t+1} - K^{n+1}_t \rangle + \frac{\bar{h}_{\text{grad}}}{2} \|K^{n+1}_{t+1} - K^{n+1}_t\|^2_F\\
    &= \langle \nabla C(K^{n+1}_t) - v^{n+1}_t,-\eta v^{n+1}_t \rangle - \eta \|v^{n+1}_t \|^2_F+ \frac{\bar{h}_{\text{grad}}}{2} \|K^{n+1}_{t+1} - K^{n+1}_t\|^2_F\\
    &\stackrel{(i)}{\leq} \frac{\eta}{2} \|\nabla C(K^{n+1}_t) - v^{n+1}_t\|^2_F - \frac{\eta}{2} \|v^{n+1}_t\|^2_F+ \frac{\bar{h}_{\text{grad}}}{2} \|K^{n+1}_{t+1} - K^{n+1}_t\|^2_F\\
    &\stackrel{(ii)}{\leq} \frac{3\eta}{4} \|\nabla C(K^{n+1}_t) - v^{n+1}_t\|^2_F - \frac{\eta}{8}\|\nabla C(K^{n+1}_t)\|^2_F - \frac{\eta}{4} \|v^{n+1}_t\|^2_F+\frac{\bar{h}_{\text{grad}}}{2} \|K^{n+1}_{t+1} - K^{n+1}_t\|^2_F\\
    &=  \frac{3\eta}{4} \|\nabla C(K^{n+1}_t) - v^{n+1}_t\|^2_F - \frac{\eta}{8}\|\nabla C(K^{n+1}_t)\|^2_F- \frac{\eta}{4} \|v^{n+1}_t\|^2_F+\frac{\bar{h}_{\text{grad}}\eta^2}{2} \|v_t^{n+1}\|^2_F\\
    &\stackrel{(iii)}{\leq}  \frac{3\eta}{4} \|\nabla C(K^{n+1}_t) - v^{n+1}_t\|^2_F - \frac{\eta}{8}\|\nabla C(K^{n+1}_t)\|^2_F- \frac{\barhgrad}{2}\|K^{n+1}_{t+1} - K^{n+1}_t\|^2_F,
\end{align*}
where $(i)$ and $(ii)$ are due to Young's inequality. For the latter we use $\|\nabla C(K^{n+1}_t)\|^2_F \leq 2\|v^{n+1}_t\|^2_F + 2 \|\nabla C(K^{n+1}_t) - v^{n+1}_t \|^2_F$. In $(iii)$ we select the step-size $\eta \leq \frac{1}{4\barhgrad}$. Thus, by now taking the expectation with respect to problem randomness we obtain
\begin{align}\label{eq:Lip_cost_difference}
    \mc{E}\left(C(K^{n+1}_{t+1}) - C(K^{n+1}_{t})\right) &\leq \frac{3\eta}{4} \mc{E}\|\nabla C(K^{n+1}_t) - v^{n+1}_t\|^2_F - \frac{\eta}{8}\mc{E}\|\nabla C(K^{n+1}_t)\|^2_F- \frac{\barhgrad}{2}\mc{E}\|K^{n+1}_{t+1} - K^{n+1}_t\|^2_F,
\end{align}
let us now control the term $\frac{3\eta}{4}\mc{E}\|\nabla C(K^{n+1}_t) - v^{n+1}_t\|^2_F$

\allowdisplaybreaks
\begin{align}\label{eq:bound_st_grad}
  &\frac{3\eta}{4}\mc{E}\|\nabla C(K^{n+1}_t) - v^{n+1}_t\|^2_F  = \frac{3\eta}{4}\mc{E}\|\nabla C(K^{n+1}_t) - \widetilde{\nabla}C(\tilde{K}^n) - \overline{\nabla} C(K^{n+1}_t) + \overline{\nabla} C(\Tilde{K}^n) \|^2_F\notag\\
  &{=}\frac{3\eta}{4}\mc{E}\|\nabla C(K^{n+1}_t) - \widetilde{\nabla}C(\tilde{K}^n)  + \frac{1}{n_2}\sum_{i=1}^{n_2} \overline{g}_i (\Tilde{K}^n) - \overline{g}_i(K^{n+1}_t)\|^2_F\notag\\
   &\stackrel{(i)}{=}\frac{3\eta}{4n^2_2}\sum_{i=1}^{n_2}\mc{E}\|\nabla C(K^{n+1}_t) - \widetilde{\nabla}C(\tilde{K}^n)  + \overline{g}_i (\Tilde{K}^n) - \overline{g}_i(K^{n+1}_t)\|^2_F\notag\\
   &=\frac{3\eta}{4n^2_2}\sum_{i=1}^{n_2}\mc{E}\|\nabla C(K^{n+1}_t) - \widetilde{\nabla}C(\tilde{K}^n) - \nabla C(\tilde{K}^n) +  \nabla C(\tilde{K}^n) + \overline{g}_i (\Tilde{K}^n) - \overline{g}_i(K^{n+1}_t)\|^2_F\notag\\
    &\stackrel{(ii)}{=}\frac{3\eta}{4n^2_2}\sum_{i=1}^{n_2}\mc{E}\|\nabla C(K^{n+1}_t) - \nabla C(\tilde{K}^n) + \overline{g}_i (\Tilde{K}^n) - \overline{g}_i(K^{n+1}_t)\|^2_F+ \frac{3\eta}{4n_2}\mc{E}\|\nabla C(\tilde{K}^n) -  \widetilde{\nabla} C(\tilde{K}^n) \|^2_F\notag\\
    &\stackrel{(iii)}{\leq}\frac{3\eta}{2n^2_2}\sum_{i=1}^{n_2}\mc{E}\|\mc{E}\left( \overline{g}_i(K^{n+1}_t) - \overline{g}_i (\Tilde{K}^n) \right) - \overline{g}_i(K^{n+1}_t) + \overline{g}_i (\Tilde{K}^n)  \|^2_F+ \frac{3\eta}{2n_2}\mc{E}\|\mc{E}\widetilde{\nabla} C(\tilde{K}^n) -  \widetilde{\nabla} C(\tilde{K}^n) \|^2_F+ \frac{15\eta\mathcal{B}(r)}{2n_2}\notag\\
    &\stackrel{(iv)}{\leq} \frac{3\eta}{2n^2_2}\sum_{i=1}^{n_2}\mc{E}\| \overline{g}_i(K^{n+1}_t) + \overline{g}_i (\Tilde{K}^n)  \|^2_F+ \frac{3\eta}{2n_2}\mc{E}\|\widetilde{\nabla} C(\tilde{K}^n) \|^2_F+ \frac{15\eta\mathcal{B}(r)}{2n_2}\notag\\
    &\stackrel{(v)}{\leq} \frac{3\eta C^2_g}{2n_2}\mc{E}\| K^{n+1}_t - \Tilde{K}^n  \|^2_F+ \frac{3\eta}{2n_2}\mc{E}\|\widetilde{\nabla} C(\tilde{K}^n) \|^2_F+ \frac{15\eta\mathcal{B}(r)}{2n_2}\notag\\
    & \stackrel{(vi)}{\leq} \frac{3\eta C^2_g}{2n_2}\mc{E}\| K^{n+1}_t - \Tilde{K}^n  \|^2_F+ \frac{3\eta}{2n_2}\left(8d^2\mathcal{B}(r)+ 2d^2\mc{E}\|\nabla C(\tilde{K}^n) \|^2_F\right)+ \frac{15\eta\mathcal{B}(r)}{2n_2}\notag \\
    & \stackrel{(vii)}{\leq} \eta \phi_1 \mc{E}\| K^{n+1}_t - \Tilde{K}^n  \|^2_F+ \frac{\eta}{16}\mc{E}\|\nabla C(K^{n+1}_t) \|^2_F +   \frac{\eta\phi_2\mathcal{B}(r)}{n_2},
\end{align}
where $(i)$ is due to the independence across the inner loop samples, since $\mc{E}\|z_1 + \ldots + z_n\|^2_F = \mc{E} \|z_1\|^2_F + \ldots +  \mc{E} \|z_n\|^2_F$ for independent and zero mean variables $\{z_1,\ldots,z_n\}$, and $(ii)$ is due to the fact that inner and outer loop samples are independent. In $(iii)$ we decompose the true gradient into estimated gradient and bias, both ZO1P and ZO2P biases are then combined together in one term. In addition, $(iv)$ is due to the fact that $\mc{E}\|X - \mc{E}X\|^2_F \leq \mc{E}\|X\|^2_F$ for any given random vector $X$. $(v)$ follows from Assumption \ref{assumption:Lipschitz_approx_gradient}. $(vi)$ is due to  Lemma \ref{lemma:two_point_estimated_grad} and for $(vii)$ we use the local Lipschitz property of the true gradient (Lemma \ref{lemma:Lipschitz_true_cost}) with the selection of $n_2 \geq 96d^2$. Here we define $\phi_1 = \frac{\frac{3C^2_g}{2} + 6\barhgrad^2d^2}{n_2}$, $\phi_2= \frac{15}{2} + 12d^2$. 

Next, \eqref{eq:bound_st_grad} can be applied in \eqref{eq:Lip_cost_difference} to obtain

\begin{align*}
    \mc{E}\left(C(K^{n+1}_{t+1}) - C(K^{n+1}_{t})\right) &\leq -\frac{\eta}{16}\mc{E}\|\nabla C(K^{n+1}_t)\|^2_F+ \frac{\eta\phi_2\mathcal{B}(r)}{n_2}+ \eta\phi_1\mc{E}\|K^{n+1}_t - \tilde{K}^n\|^2_F
    - \frac{\barhgrad}{2}\mc{E}\|K^{n+1}_{t+1} - K^{n+1}_t\|^2_F,
\end{align*}

We can now use the gradient domination of the cost function (Lemma \ref{lemma:PL_condition}) to write
\begin{align*}
    \mc{E}\left(C(K^{n+1}_{t+1}) - C(K^*)\right) &\leq \gamma\mc{E}\left(C(K^{n+1}_{t}) - C(K^*)\right)+ \frac{\eta\phi_2\mathcal{B}(r)}{n_2}+ \eta\phi_1\mc{E}\|K^{n+1}_t - \tilde{K}^n\|^2_F
    - \frac{\barhgrad}{2}\mc{E}\|K^{n+1}_{t+1} - K^{n+1}_t\|^2_F,
\end{align*}
with $\gamma := 1 - \frac{\eta \lambda}{16}$. Observe that by applying Young's inequality we can write
\begin{align*}
    \|K^{n+1}_{t+1} - K^{n+1}_t\|^2_F &\geq \frac{1}{1+\alpha}\|K_{t+1}^{n+1} - \tilde{K}^n\|^2_F- \frac{1}{\alpha}\|K^{n+1}_t - \tilde{K}^n\|^2_F,
\end{align*}
for any $\alpha > 0$, which implies
\begin{align}\label{eq:telecoping_ref}
    \mc{E}\left(C(K^{n+1}_{t+1}) - C(K^*)\right) &\leq \gamma\mc{E}\left(C(K^{n+1}_{t}) - C(K^*)\right)+ \left(\eta\phi_1 + \frac{\barhgrad}{2\alpha}\right)\mc{E}\|K^{n+1}_t - \tilde{K}^n\|^2_F 
    - \frac{\barhgrad}{2(1+\alpha)}\mc{E}\|K^{n+1}_{t+1} - \tilde{K}^n\|^2_F\notag\\
    & + \frac{\eta\phi_2\mathcal{B}(r)}{n_2},
\end{align}
which leads to 
\begin{align*}
    \mc{E}\left(C(K^{n+1}_{T}) - C(K^*)\right) &\leq \gamma^T\mc{E}\left(C(K^{n+1}_{0}) - C(K^*)\right) + \frac{\eta\phi_2\mathcal{B}(r)}{n_2}\sum_{t=0}^{T-1}\gamma^t +\mathcal{C}\\
     &= \gamma^T\mc{E}\left(C(K^{n+1}_{0}) - C(K^*)\right) 
     + \frac{16\phi_2\mathcal{B}(r)(1-\gamma^T)}{n_2\lambda} + \mathcal{C},
\end{align*}
with 
\begin{align*}
    \mathcal{C} &= \sum_{t=0}^{T-1}\gamma^
    {T-1-t}\left[\left(\eta\phi_1 + \frac{\barhgrad}{2\alpha}\right)\mc{E}\|K^{n+1}_t - \tilde{K}^n\|^2_F-  \frac{\barhgrad\mc{E}\|K^{n+1}_{t+1} - \tilde{K}^n\|^2_F}{2(1+\alpha)}\right],
\end{align*}
by telescoping \eqref{eq:telecoping_ref} over $t=\{0,1,\ldots,T-1\}$. Let us now proceed by controlling $\mathcal{C}$. To achieve this, we first select $\alpha = 2t+1$ to write
\begin{align*}
    \mathcal{C} &\stackrel{(i)}{\leq} \sum_{t=1}^{T-1}\left(\eta\phi_1 + \frac{\barhgrad}{2(2t+1)}\right)\mc{E}\|K^{n+1}_t - \tilde{K}^n\|^2_F- \sum_{t=0}^{T-2} \frac{\barhgrad \mc{E}\|K^{n+1}_{t+1} - \tilde{K}^n\|^2_F}{4(t+1)} - \sum_{t=0}^{T-2} \frac{\barhgrad \mc{E}\|K^{n+1}_{t+1} - \tilde{K}^n\|^2_F}{4T}\\
    &\leq\sum_{t=1}^{T-1}\left(\eta\phi_1 + \frac{\barhgrad}{2(2t+1)}\right)\mc{E}\|K^{n+1}_t - \tilde{K}^n\|^2_F- \sum_{t=0}^{T-2} \frac{\barhgrad \mc{E}\|K^{n+1}_{t+1} - \tilde{K}^n\|^2_F}{4(t+1)}\\
    &\stackrel{(ii)}{\leq} \sum_{t=1}^{T-1}\left(\eta\phi_1 + \frac{\barhgrad}{2(2t+1)}\right)\mc{E}\|K^{n+1}_t - \tilde{K}^n\|^2_F- \sum_{t=1}^{T-1} \frac{\barhgrad \mc{E}\|K^{n+1}_{t} - \tilde{K}^n\|^2_F}{4t}\\
    &= \sum_{t=1}^{T-1}\left(\eta\phi_1 + \frac{\barhgrad}{2(2t+1)} 
    - \frac{\barhgrad}{4 t}\right)\mc{E}\|K^{n+1}_t - \tilde{K}^n\|^2_F\\
    &= \sum_{t=1}^{T-1}\left(\eta\phi_1 - \frac{\barhgrad}{4(2t+1)t}\right)\mc{E}\|K^{n+1}_t - \tilde{K}^n\|^2_F,
\end{align*}
where $(i)$ is due to the fact that $K^{n+1}_0 = \tilde{K}^n$ and $0 < \gamma < 1$. $(ii)$ follows from a simple change of variables $\bar{t}=t+1$. Now, by carefully selecting the inner-loop samples $n_2$ according to $n_2 \geq \frac{\left(3C^2_g + 12\barhgrad^2d^2\right)T^2}{\barhgrad^2}$ we obtain $\mathcal{C}\leq 0$, which implies
\begin{align*}
    \mc{E}\left(C(K^{n+1}_{T}) - C(K^*)\right) &\leq  \gamma^T\mc{E}\left(C(K^{n+1}_{0}) - C(K^*)\right) + \frac{16\phi_2\mathcal{B}(r)(1-\gamma^T)}{n_2\lambda}, 
\end{align*}
which can be unrolled over $n=\{0,\ldots,N-1\}$ to obtain
\begin{align*}
    \mc{E}\left(C(K_{\text{out}}) - C(K^*)\right) &\leq  \gamma^{NT}\mc{E}\left(C(K_{0}) - C(K^*)\right)+ \frac{\mathcal{B}(r) \left(120+ 192d^2\right)}{n_2\lambda},
\end{align*}
which completes the proof.

\subsection{Numerical Experiments - Additional Details} \label{appendix:exp_additional_details}

In our numerical experiments we consider a unstable discrete-time and LTI system as in \eqref{eq:LTI_system} with system and cost matrices described by

\begin{align*}
    A = \begin{bmatrix}
        1.20 & 0.50 & 0.40\\
        0.01 & 0.75 & 0.30\\
        0.10 & 0.02 & 1.50
    \end{bmatrix},  B = \begin{bmatrix}
        \frac{1}{2} \\ 1 \\ \frac{1}{2}
    \end{bmatrix},  Q=2I_3,  R = \frac{1}{2}.
\end{align*}

In addition, we initialize Algorithm \ref{algorithm:LQR_ZO2P} and \ref{algorithm:LQR_SVRPG} with the sub-optimal stabilizing controller $K_0 = [0.15 ;\ -0.45 ;\ 3.80]$.

\end{document}